\documentclass[12pt]{amsart}
\usepackage{amscd,amssymb,epsf}
\usepackage[mathscr]{eucal}

\newtheorem*{thm*}{Theorem}

\newtheorem*{corollary*}{Corollary}
\newtheorem*{claim*}{Claim}

\numberwithin{equation}{subsection}
\numberwithin{thm}{subsection}
\begin{document}
\newcommand{\R}{{\mathbb R}}
\newcommand{\C}{{\mathbb C}}
\newcommand{\Z}{{\mathbb Z}}
\newcommand{\B}{{\mathbb B}}
\newcommand{\A}{{\mathbb A}}
\renewcommand{\aa}{{\mathfrak a}}
\newcommand{\Id}{{\mathbb I}}
\renewcommand{\P}{{\mathbb P}}
\newcommand{\cpo}{{\C\P^1}}
\newcommand{\Hom}{\o{Hom}}
\newcommand{\Tei}{{\mathfrak T}}
\newcommand{\Teich}{\Tei_M}
\newcommand{\Hh}{{\mathfrak H}}
\newcommand{\Cc}{{\mathfrak C}}
\newcommand{\Hm}[1]{\Hom (\pi ,#1)}
\newcommand{\Hmm}[1]{\Hom (\pi ,#1)/#1}
\newcommand{\Hmmm}[1]{\Hom (\pi ,#1)/\hspace{-3pt}/#1}
\newcommand{\hmg}{\Hom (\pi ,G)}
\renewcommand{\k}{\mathbf{k}}
\newcommand{\Ker}{\o{Ker}}
\newcommand{\Aut}{\o{Aut}}
\newcommand{\Ad}{\o{Ad}}
\newcommand{\hg}{\Hmmm{G}}
\newcommand{\tr}{{\bf\mathsf{tr}}}
\newcommand{\dd}[1]{\frac{\partial}{\partial{#1}}}
\newcommand{\SL}[1]{{\o{SL}}({#1})}
\newcommand{\SLt}{{\SL{2}}}
\newcommand{\GL}[1]{{\o{GL}}({#1})}
\newcommand{\slt}{{\SL{2,\C}}}
\newcommand{\slr}{{\SL{2,\R}}}
\newcommand{\sothc}{{\o{SO}{(3,\C)}}}
\newcommand{\soth}{{\o{SO}(3)}}
\newcommand{\soto}{{\o{SO}(2,1)}}
\newcommand{\sot}{{\o{SO}{(2)}}}
\newcommand{\sooo}{{\o{SO}{(1,1)}}}
\newcommand{\Ht}{\mathbf{H}^2}
\newcommand{\slthr}{{\mathfrak{sl}(3,\R)}}
\newcommand{\rpt}{{\R\P}^2}
\newcommand{\PSLtC}{{\mathsf{PSL}(2,\C)}}
\newcommand{\SLthR}{{\mathsf{SL}(3,\R)}}
\newcommand{\SLthZ}{{\mathsf{SL}(3,\Z)}}
\newcommand{\PGLthR}{{\mathsf{PGL}(3,\R)}}

\renewcommand{\ll}{{\mathcal L}}

\title{Bulging deformations of convex $\rpt$-manifolds}
\author{William M.~Goldman}
\address{ Mathematics Department,
University of Maryland, College Park, MD  20742 USA  }
\email{ wmg@math.umd.edu }

\thanks{The author gratefully acknowledges partial support from
National Science Foundation grant }
\date{\today}
\subjclass{57M05 (Low-dimensional topology), 20H10 (Fuchsian groups and their
generalizations), 30F60 (Teichm\"uller theory)}
\keywords{
surface, $\rpt$-structure, measured lamination, earthquake, 
Teichm\"uller space}
\begin{abstract}
We define deformations of convex $\rpt$-surfaces.
\end{abstract}
\maketitle A {\em convex $\rpt$-manifold\/} is a representation of a
surface $S$ as a quotient $\Omega/\Gamma$, where $\Omega\subset\rpt$
is a convex domain and $\Gamma\subset\SLthR$ is a discrete group of
collineations acting properly on $\Omega$. We shall describe a
construction of deformations of such structures based on Thurston's
earthquake deformations for hyperbolic surfaces and {\em quakebend
deformations \/} for $\cpo$-manifolds.

\begin{figure} [htp]
\centerline{\epsfxsize=3in \epsfbox{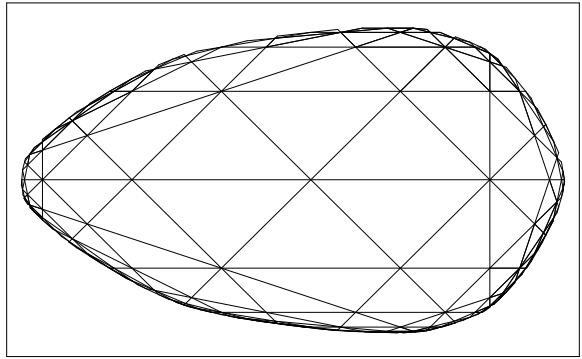}}
\caption{A convex domain tiled by triangles}
\label{fig:t334}
\end{figure}

In general if $\Omega/\Gamma$ is a convex $\rpt$-manifold which is a
{\em closed\/} surface $S$ with $\chi(S)$, then 
either $\partial\Omega$ is a conic, or
$\partial\Omega$ is a $C^1$ convex curve (Benz\'ecri~\cite{Benzecri})
which is not $C^2$ (Kuiper~\cite{Kuiper}). 
In fact its derivative is H\"older continuous with H\"older exponent strictly
between 1 and 2. Figure~\ref{fig:t334} depicts such a domain tiled
by the $(3,3,4)$-triangle tesselation.

\newpage
This drawing actually arises from Lie algebras (see
Kac-Vinberg~\cite{KacVinberg}). Namely the Cartan matrix
\begin{equation*}
C = \bmatrix 2 & -1 & -1 \\ -2 & 2 & -1 \\ -1 & -1 & 2 \endbmatrix    
\end{equation*}
determines a group of reflections as follows. For $i=1,2,3$ let
$E_{ii}$ denote the elementary matrix having entry $1$ in the $i$-th
diagonal slot.  Then, for $i=1,2,3$, the reflections 
\begin{equation*}
\rho_i = I - E_{ii} C  
\end{equation*}
generate a discrete subgroup of $\SLthZ$ which acts
properly on the convex domain depicted in Figure~\ref{fig:t334}
(and appears on the cover of the November 2002 Notices of the American
Mathematical Society).. This group is the Weyl group of a 
hyperbolic Kac-Moody Lie algebra.

We describe here a general construction of such convex domains as limits
of {\em piecewise conic\/} curves.

If $\Omega/\Gamma$ is a convex $\rpt$-manifold homeomorphic to a closed
esurface $S$ with $\chi(S)<0$, then every element $\gamma\in\Gamma$ is
{\em positive hyperbolic,\/} that is, conjugate in $\SLthR$ to a diagonal
matrix of the form
\begin{equation*}
\delta = \bmatrix  e^s & 0 & 0 \\ 0 & e^t & 0 \\ 0 & 0 & e^{-s-t} \endbmatrix.
\end{equation*}
where $ s > t > -s -t $. Its centralizer is the {\em maximal $\R$-split torus\/}
$\A$ consisting of all diagonal matrices in $\SLthR$. It is isomorphic
to a Cartesian product $\R^*\times \R^*$ and has four connected components.
Its identity component $\A^+$ 
consists of diagonal matrices with positive entries.

The {\em roots\/} are linear functionals on its Lie algebra $\aa$, the
{\em Cartan subalgebra.\/} Namely, $\aa$ consists of diagonal matrices
\begin{equation}\label{eq:diagonalmatrix}
a  = \bmatrix  a_1 & 0 & 0 \\ 0 & a_2 & 0 \\ 0 & 0 & a_3 \endbmatrix. 
\end{equation}
where $a_1 + a_2 + a_3=0$. The roots are the six linear functionals
on $\aa$ defined by
\begin{equation*}
a \stackrel{\alpha_{ij}}\longmapsto  a_i - a_j 
\end{equation*}
where $1\le i \neq j \le 3$. Evidently $\alpha_{ji}=-\alpha_{ij}$.

Writing
$a(s,t)$ for the diagonal matrix \eqref{eq:diagonalmatrix} with
\begin{equation*}
a_1 = s, \qquad a_1 = 2, \qquad  a_3 = - s - t,
\end{equation*}
the roots are the linear functionals defined by 
\begin{align*}
a(s,t) & \stackrel{\alpha_{12}}\longmapsto  s -t \\
a(s,t) & \stackrel{\alpha_{21}}\longmapsto  t -s \\
a(s,t) & \stackrel{\alpha_{23}}\longmapsto  t - (- s -t) = s + 2 t \\
a(s,t) & \stackrel{\alpha_{32}}\longmapsto  (- s -t) - t = - s - 2 t \\
a(s,t) & \stackrel{\alpha_{31}}\longmapsto  (- s -t) - s  = -2 s -t \\
a(s,t) & \stackrel{\alpha_{13}}\longmapsto  s - (- s -t)  = 2 s  + t 
\end{align*}
which we write as
\begin{align*}
\alpha_{12} &= \bmatrix 1 &  -1 \endbmatrix \\
\alpha_{21} &= \bmatrix -1 &  1 \endbmatrix \\
\alpha_{23} &= \bmatrix 1 &  2 \endbmatrix \\
\alpha_{32} &= \bmatrix -1 & -2 \endbmatrix \\
\alpha_{31} &= \bmatrix -2 & -1 \endbmatrix \\
\alpha_{13} &= \bmatrix 2 & 1 \endbmatrix 
\end{align*}

The {\em Weyl group\/} is generated by reflections in the roots and
in this case is just the symmetric group, consisting of permutations
of the three variables $a_1,a_2,a_3$ in $a$ (as in \eqref{eq:diagonalmatrix}).
A fundamental domain is the {\em Weyl chamber\/} consisting of
all $a$ satisfying $\alpha_{12}>0$ and $\alpha_{23}>0$.
This corresponds to the ordering of the roots where 
$\alpha_{12}> \alpha_{23}$ are the {\em positive simple roots.\/} 
In other words,  the roots are totally ordered by:
the rule 
\begin{equation*}
\alpha_{13} >  \alpha_{12}> \alpha_{23} > 0 >  \alpha_{32} > \alpha_{21} 
> \alpha_{21} > \alpha_{31}.
\end{equation*}
In terms of the parametrization of $\aa$ by  $a(s,t)$, the Weyl chamber
equals
\begin{equation*}
\{ a(s,t) \mid s \ge t \ge -\frac12 s \}.
\end{equation*}

The {\em trace form\/} on $\slthr$ defines the inner product 
$\langle,\rangle$ with associated quadratic form
\begin{equation*}
\tr \big(a(s,t)^2\big) = 2 (s^2 + s t + t^2) = 2 \vert s + \omega t\vert^2
\end{equation*}
where $\omega = \frac12 + \frac{\sqrt{-3}}2 = e^{\pi i/3}$ 
is the primitive sixth root of $1$.

The elements of $\slthr$ which dual to the roots (via the inner product 
$\langle,\rangle$) are the {\em root vectors:\/}
\begin{align*}
h_{12} &= a(1,-1),\\
h_{21} &= a(-1,1), \\
h_{23} &= a(0,1),\\
h_{32} & = a(0,-1),\\
h_{31} & = a(-1,0),\\
h_{13} &= a(1,0) 
\end{align*}
The Weyl chamber consists of all
\begin{equation*}
a(s,\lambda s) = \bmatrix s & 0 & 0 \\ 0 & \lambda s & 0 \\ 0 & 0 
& -(1+\lambda) s \endbmatrix
\end{equation*}
where $1 \ge \lambda \ge -\frac12$. Its boundary consists of the rays
generated by the {\em singular elements\/}
\begin{equation*}
a(1,1) = h_{13} + h_{23} = h_{12} + 2  h_{23}
\end{equation*}
and
\begin{equation*}
a(2,-1) = h_{12} + h_{13} = 2  h_{12} +  h_{23}.
\end{equation*}
The sum of the simple positive roots is the element
\begin{equation*}
a(1,0) = h_{13}= h_{12} + h_{23}
\end{equation*}
which generates the one-parameter subgroup
\begin{equation*}
H_t := \exp\big(a(t,0)\big)  = \bmatrix e^t & 0 & 0 \\ 0 & 1 & 0 \\ 0 & 0 & 
e^{-t} \endbmatrix.
\end{equation*}

\begin{figure}[htp]
\centerline{\epsfxsize=3in \epsfbox{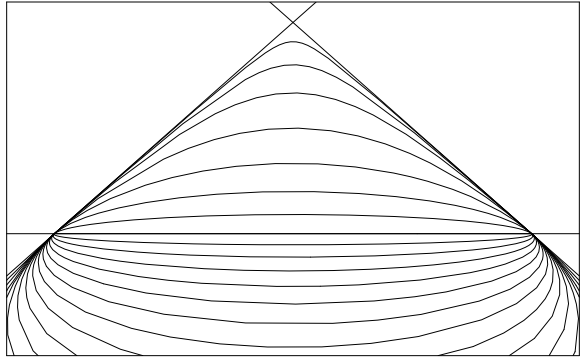}}
\caption{Conics tangent to a triangle}
\label{fig:conicstriangle}
\end{figure}

The orbits of $\A^+$ on $\rpt$ are the four open 2-simplices defined by
the homogeneous coordinates, their (six) edges and their (three) vertices.
The orbits of $H_t$ are arcs of conics depicted in 
Figure~\ref{fig:conicstriangle}.

Associated to any measured geodesic lamination $\lambda$ on a hyperbolic
surface $S$ is {\em bulging deformation\/} as an $\rpt$-surface.
Namely, one applies a one-parameter group of collineations
\begin{equation*}
\bmatrix 1 & 0 & 0 \\ 0 & e^t & 0 \\ 0 & 0 & 1 \endbmatrix
\end{equation*}
to the coordinates on either side of a leaf.
This extends Thurston's {\em earthquake\/} deformations 
(the analog of {\em Fenchel-Nielsen twist deformations\/} along
possibly infinite geodesic laminations), and the {\em bending deformations \/}
in $\PSLtC$.

In general, if $S$ is a convex $\rpt$-manifold, then deformations are
determined by a geodesic lamination with a transverse measure taking values
in the Weyl chamber of $\SLthR$. When $S$ is itself a hyperbolic surface,
all the deformations in the singular directions become earthquakes
and deform $\partial\tilde S$ trivially (just as in $\PSLtC$.

\begin{figure}[htp]
\centerline{\epsfxsize=3in \epsfbox{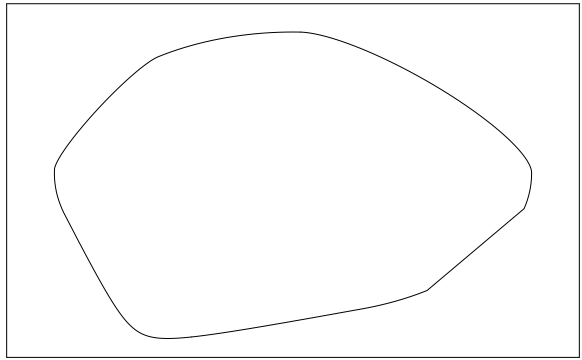}}
\caption{Deforming a conic}
\label{fig:bulgedconic}
\end{figure}

\newpage
\begin{figure}[htp]
\centerline{\epsfxsize=3in \epsfbox{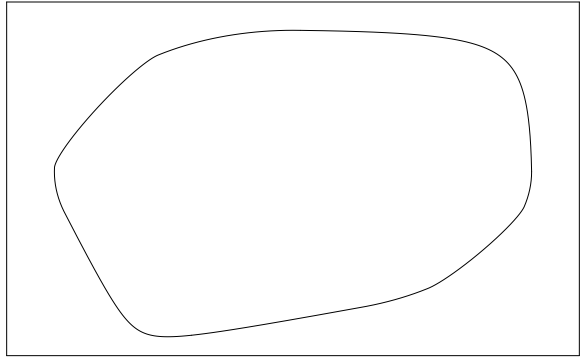}}
\caption{A piecewise conic}
\label{A piecewise conic}
\end{figure}

\bigskip
\begin{figure}[htp]
\centerline{\epsfxsize=3in \epsfbox{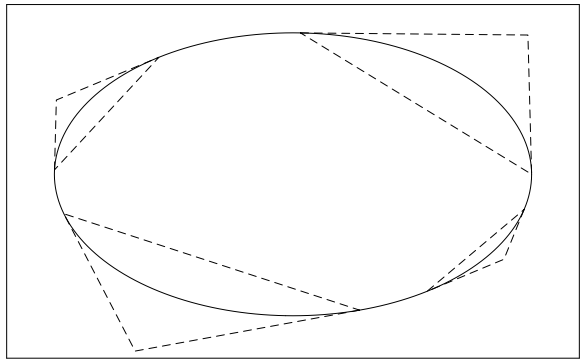}}
\caption{Bulging data}
\label{Bulging data}
\end{figure}

\bigskip
\begin{figure}[htp]
\centerline{\epsfxsize=3in \epsfbox{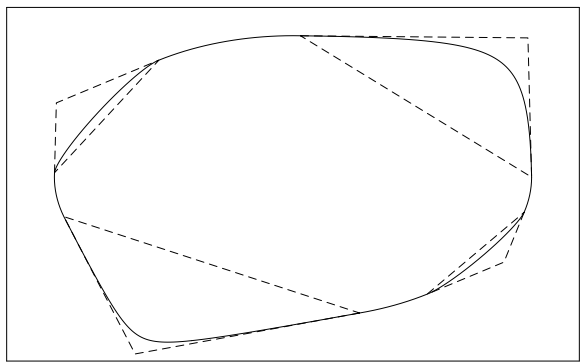}}
\caption{The deformed conic}
\label{The deformed conic}
\end{figure}

\bigskip
\begin{figure}[htp]
\centerline{\epsfxsize=3in \epsfbox{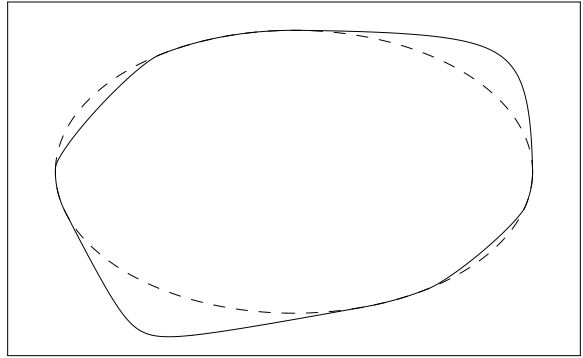}}
\caption{The conic with its deformation}
\label{The conic with its deformation}
\end{figure}

\makeatletter \renewcommand{\@biblabel}[1]{\hfill#1.}\makeatother

\end{document}